\documentclass{article}
\usepackage{amsfonts}
\usepackage{amsmath}
\usepackage[all]{xy}
\newtheorem{theorem}{Theorem}[section]
\newtheorem{definition}[theorem]{Definition}

\newtheorem{corollary}[theorem]{Corollary}

\newtheorem{lemma}[theorem]{Lemma}

\usepackage{mathrsfs}
\usepackage{amssymb}
\usepackage{amscd}
\date{}

\usepackage{amssymb}
\usepackage{extarrows}
\usepackage{cite}

\begin{document}

\title{Gr\"{o}bner-Shirshov bases for Lie $\Omega$-algebras and free Rota-Baxter Lie algebras\footnote{Supported by the NNSF of China (11171118, 11571121).}}

\author{
Jianjun Qiu and Yuqun Chen\footnote {Corresponding author.}   \\
{\small \ School of Mathematical Sciences, South China Normal
University}\\
{\small Guangzhou 510631, P. R. China}\\
{\small jianjunqiu@126.com}\\
{\small yqchen@scnu.edu.cn}
}

\maketitle \noindent\textbf{Abstract:} In this paper, we generalize the  Lyndon-Shirshov
words to Lyndon-Shirshov $\Omega$-words on a set $X$ and prove that the set of all non-associative  Lyndon-Shirshov $\Omega$-words forms a linear basis of the free  Lie $\Omega$-algebra on the set $X$. From this, we establish   Gr\"{o}bner-Shirshov bases theory for Lie $\Omega$-algebras.  As applications, we give Gr\"{o}bner-Shirshov bases for   free   $\lambda$-Rota-Baxter Lie algebras, free  modified $\lambda$-Rota-Baxter Lie algebras and free Nijenhuis Lie algebras and then  linear bases of such three free  algebras  are obtained.

 \ \


  \ \


\section{Introduction}

  Let $(R, \cdot)$ be an algebra over a field  $k$, $\lambda\in k$ and    $ P:R\rightarrow R$  a linear map satisfying
$$
 P(x)\cdot P(y)  = P(  P(x)\cdot y)  +P( x\cdot P(y))+\lambda P^n ( x\cdot y ),  \    x,y \in  R,
$$
where $n$ is a nonnegative integer.   If $n=1$ (resp, $n=0$, $n=2$ and $\lambda=-1$), then the operator $P$ is called a  $\lambda$-Rota-Baxter (resp. modified $\lambda$-Rota-Baxter, Nijenhuis) operator on the algebra $(R, \cdot)$ and
the triple $(R, \cdot,  P)$ is called a    $\lambda$-Rota-Baxter  (resp. modified $\lambda$-Rota-Baxter, Nijenhuis)  algebra.

Rota-Baxter operator on an associative algebra  was introduced  by  G.  Baxter    to solve an analytic  problem in probability  \cite{Bax60} and    was studied by  G.-C. Rota  \cite{Ro69} in combinatorics later.
The  Rota-Baxter operator of weight $\lambda=0$ on  Lie algebra  is also   called the operator form of the classical Yang-Baxter equation due to Semenov-Tian-Shansky's work \cite{sts}.  Modified Rota-Baxter operator  on  an associative algebra was introduced  by  K. Ebrahimi-Fard \cite{EF02} with motivation from
modified classical Yang-Baxter equation on a  Lie algebra \cite{sts}.  The Nijenhuis operator on an associative algebra was introduced by J. Cari\~{n}ena et al. \cite{car00}
to study quantum bi-Hamiltonian systems. In \cite{uc10}, Nijenhuis operators are constructed by analogy
with Poisson-Nijenhuis geometry, from relative Rota-Baxter operators. Nijenhuis operators on Lie algebras play an important role in the study
of integrability of nonlinear evolution equations \cite{dor93}.
Recently,   there are  some results on  Rota-Baxter operators on Lie algebras and related topic, for example,   see \cite{amm, ab, bgn10, lhb, pgb}.

There are many constructions of free  $\lambda$-Rota-Baxter associative algebras, free modified $\lambda$-Rota-Baxter associative algebras and free Nijenhuis associative algebras   by using  different methods, for example, \cite{AM,  bcq, Ca72, EF04,  EG08a, EG08b, gg15,  lguo, Gu09,  GK00a, GK00b,   GK08, lg12, Ro69}.  However, as we know,  there are no constructions of  free $\lambda$-Rota-Baxter Lie algebras, free  modified $\lambda$-Rota-Baxter Lie algebras and free Nijenhuis Lie  algebras. We will  apply  Gr\"{o}bner-Shirshov bases method to construct such three free   Lie algebras.

Gr\"{o}bner bases and Gr\"{o}bner-Shirshov bases  were invented
independently by A.I. Shirshov for ideals of free (commutative,
anti-commutative) non-associative algebras \cite{Sh62a,Shir3}, free
Lie algebras \cite{Sh,Shir3} and implicitly free associative
algebras \cite{Sh,Shir3}  (see also \cite{Be78,Bo76}), by H.
Hironaka \cite{Hi64} for ideals of the power series algebras (both
formal and convergent), and by B. Buchberger \cite{Bu70} for ideals
of the polynomial algebras. Gr\"{o}bner bases and
Gr\"{o}bner-Shirshov bases theories have been proved to be very
useful in different branches of mathematics, including commutative
algebra and combinatorial algebra. It is a powerful tool to solve
the following classical problems: normal form; word problem;
conjugacy problem; rewriting system; automaton; embedding theorem;
PBW theorem;  extension; homology;  growth function; Dehn function;
complexity; etc. See, for example, the books \cite{AL, BKu94, BuCL,
BuW, CLO, Ei}, the papers \cite{Be78, b72,Bo76,BCM, DK10, gg15,  Mikhalev92}, and the
surveys \cite{BC14, BCS, BFKK00, BK03}.

The concept of $\Omega$-algebra   was
introduced  by A.G. Kurosh \cite{Ku60} under an influence of the
concept of multioperator group of  P.J. Higgins \cite{Hi56}. An  $\Omega$-algebra  over a field $k$   is a $k$-algebra $A$  with  a set of  multilinear operators $\Omega$, where $\Omega=\bigcup_{m=1}^{\infty}\Omega_{m}$
and each $\Omega_{m}$ is a set of $m$-ary multilinear operators on $A$.  $\lambda$-Rota-Baxter Lie algebras, modified $\lambda$-Rota-Baxter Lie algebras and  Nijenhuis Lie  algebras are Lie $\Omega$-algebras with one operator. In  \cite{bcq},
the associative $\Omega$-words on a set $X$ are introduced, it is shown that the set of all associative $\Omega$-words on  $X$ forms a linear basis of the
free  associative $\Omega$-algebra on  $X$, and then Gr\"{o}bner-Shirshov bases theory  for associative $\Omega$-algebras is established.
For    Gr\"{o}bner-Shirshov bases  theory of  various $\Omega$-algebras and their applications,     see  \cite{ bch, bcq, qc, DH08,  gg15}.

The first linear basis of a free Lie algebra $Lie(X)$ on a set $X$ had been given
by M. Hall 1950 \cite{M.H50}. K.-T. Chen, R.H. Fox, R.C. Lyndon
  1958 \cite{CFL58} and
A.I. Shirshov 1958 \cite{Sh} introduced non-associative Lyndon-Shirshov
words on  $X$ and proved that they form a linear basis of $Lie(X)$, independently.  One of the main
applications of Lyndon-Shirshov basis  is the Shirshov's theory of
Gr\"{o}bner-Shirshov bases theory for Lie algebras \cite{Sh}.

The paper is organized as follows. In section 2, we review the concept and some related properties of Lyndon-Shirshov words on a set $X$, generalize associative (resp. non-associative)  Lyndon-Shirshov words to associative (resp. non-associative)  Lyndon-Shirshov $\Omega$-words on a set $X$ and show that the set of all non-associative  Lyndon-Shirshov $\Omega$-words forms a linear basis of the free  Lie $\Omega$-algebra on the set $X$.  In section 3, we review Gr\"{o}bner-Shirshov
bases theory   for    associative  $\Omega$-algebras, give definition of Gr\"{o}bner-Shirshov
bases   for    Lie  $\Omega$-algebras and establish Composition-Diamond lemma for  Lie $\Omega$-algebras. In section 4, we give   Gr\"{o}bner-Shirshov bases for  free $\lambda$-Rota-Baxter Lie algebras, free  modified $\lambda$-Rota-Baxter Lie algebras,  free Nijenhuis Lie  algebras  and then  linear bases of such three free algebras are obtained by Composition-Diamond lemma for  Lie $\Omega$-algebras.

\section{Free  Lie $\Omega$-algebras}

\subsection{Lyndon-Shirshov words}

In this subsection, we review the concept and some properties of Lyndon-Shirshov words, which can be found in \cite{bc07, Sh, Shir3}.

For any set $X$, we define the following notations:

$S(X)$: the  set    of all    nonempty  associative words  on  $X$.

$X^*$:  the  set    of all   associative words  on  $X$ including the empty word $1$.

$X^{**}$: the  set    of all  non-associative words  on  $X$.

For any $u\in X^*$, denote $deg(u)$ to be  the degree (length) of $u$.  Let $>$  be  a   well order on $X$.    Define the lex-order $>_{lex}$ and the deg-lex order $>_{deg-lex}$ on $X^*$ with respect to $>$ by:

 (i) $1>_{lex} u$ for  any nonempty word  $u$, and  if $u=x_iu'$ and $v=x_jv'$, where $x_i, x_j\in X$,   then $u>_{lex} v$ if   $x_i> x_j$, or $x_i=x_j$ and $u'>_{lex} v'$ by induction.

 (ii) $u>_{deg-lex}v$ if $deg(u)>deg(v)$, or $deg(u)=deg(v)$ and $u>_{lex}v$.

A nonempty associative  word $w$ is called  an associative Lyndon-Shirshov word on $X$,  if $  w=uv >_{lex} vu$ for any  decomposition of $w=uv$, where $1\neq  u, v\in X^*$.

A non-associative word $(u)\in X^{**}$ is said to be   a non-associative Lyndon-Shirshov word on $X$ with respect to the  lex-order $>_{lex}$, denoted by $[u]$,   if
  \begin{itemize}
    \item [(a)] $u$ is an associative Lyndon-Shirshov word on $X$;
    \item [(b)] if $(u)=((v)(w))$, then both $(v)$ and $(w)$ are non-associative Lyndon-Shirshov words on $X$;
    \item [(c)] if $(v)=((v_1)(v_2))$, then $v_2 \leq_{lex} w$.
  \end{itemize}

Denote the set of all associative (resp. non-associative) Lyndon-Shirshov words on $X$ with respect to the lex-order $>_{lex}$  by $ALSW(X)$    (resp. $NLSW(X))$.
It is well known that for any  $u\in ALSW(X)$,   there exists a unique Shirshov standard  bracketing way  $[u]$ such that $[u] \in NLSW(X)$.
Then
$
NLSW(X )=\{[u] | u\in ALSW(X )\}.
$


Let $k\langle X\rangle$  be the free associative algebra on $X$ over a field $k$ and  $Lie(X)$ be  the Lie subalgebra of   $k\langle X\rangle$  generated by   $X$ under the Lie bracket $(u v)=uv-vu$. It is well known that $Lie(X)$ is a free Lie algebra on the set $X$ with a linear basis $NLSW(X)$.

For any  $f\in k\langle X\rangle$, let $\bar{f}$ be the  leading word of $f$ with respect to the deg-lex  order $>_{deg-lex}$ on $X^*$.

\begin{lemma}\label{le2.1} (\cite{bc07, Sh, Shir3})
For any non-associative word $(u)\in X^{**}$, $(u)$ has a representation
$$
(u)=\sum \alpha_i[u_i],
$$
in $Lie(X)$, where each $\alpha_i\in k$ and  $u_i\in ALSW(X)$.
\end{lemma}

\begin{lemma}\label{le2.2}(\cite{bc07, Sh, Shir3})
If  $u\in  ALSW(X)$, then $\overline{[u]}=u$.
\end{lemma}

\begin{lemma}(\cite{bc07, Sh, Shir3})\label{lem2.3}
For any $u\in X^*$, there exists  a unique decomposition
$$
u=u_1u_2\cdots u_m,
$$
where each $u_i\in ALSW(X)$ and $u_j\leq_{lex}u_{j+1}$, $1 \leq j\leq m-1$.
\end{lemma}

\begin{lemma} \label{le2.5}(\cite{bc07,  Sh,   Shir3})  Let    $u=avb$,  where $u, v \in ALSW(X)$ and  $a, b \in X^*$. Then
$$
[u]=[a[vc]d],
$$
where $b=cd$ for some $c,d\in X^*$. Let
$$
[u]_{v}=[u]_{[vc]\mapsto [\cdots[[[v][c_1]][c_2]]\cdots [c_m]]}
$$
where $c=c_1c_2\cdots c_m$ with each  $c_i\in ALSW(X)$ and $c_{i} \leq_{lex} c_{i+1}$.  Then,
$$
[u]_v=a[v]b+\sum \alpha_ia_i[v]b_i,
$$
where each $\alpha_i\in k, a_i, b_i\in X^*$ and $a_ivb_i<_{deg-lex} avb=u$. It follows that   $\overline{[u]_{v}}=u$.
\end{lemma}

\subsection{Lyndon-Shirshov  $\Omega$-words}

In this subsection, we   define the Lyndon-Shirshov  $\Omega$-words on a set $X$.

Let
$$
\Omega=\bigcup_{m=1}^{\infty}\Omega_{m},
$$
where $\Omega_{m}$ is a set of $m$-ary operators for any $m\geq 1$. For any set
 $Y$, denote
$$
\Omega(Y):=\bigcup_{m=1}^{\infty}\left\{ \omega^{(m)}(y_1, y_2, \cdots, y_m)|y_i\in Y,  1 \leq i \leq m,  \omega^{(m)}\in \Omega_{m} \right\}.
$$

Let $X$ be a set. Define
$
\langle \Omega;  X \rangle_0:=S(X)$ and $(\Omega;  X)_0:=X^{**}$. Assume that we have defined $\langle \Omega;  X \rangle_{n-1}$ and $(\Omega;  X)_{n-1}$. Define
\begin{eqnarray*}
\langle \Omega;  X \rangle_n&:=&S(X\cup \Omega(\langle \Omega;  X \rangle_{n-1})),\\
(\Omega;  X)_n&:=& (X\cup \Omega((\Omega;  X)_{n-1}))^{**}.
\end{eqnarray*}
Then it is easy to see that for any $n\geq 0$, 
$$
\langle \Omega;  X \rangle_{n }\subseteq \langle \Omega;  X \rangle_{n+1}, \ \ (\Omega;  X)_{n } \subseteq (\Omega;  X)_{n+1}.  
$$
Denote
$$\langle \Omega;  X \rangle: = \bigcup_{n=0}^{\infty}\langle \Omega;  X \rangle_n, \ \ \ 
 (\Omega;  X): = \bigcup_{n=0}^{\infty}(\Omega;  X)_n.
$$
The elements of $\langle \Omega;  X \rangle$ (resp. $(\Omega;  X)$) are called the associative  (resp. non-associative) $\Omega$-words on $X$.

If   $u \in X \cup \Omega( \langle \Omega;  X \rangle)$, then $u$  is called   prime. Therefore, for any $u\in \langle \Omega;  X \rangle$, $u$ can be expressed   uniquely in the canonical form
$$
u=u_1u_2\cdots u_n,\ n\geq 1,
$$
where each $u_i$ is prime. The number  $n$  is called the breath of $u$, which is    denoted  by $bre(u)$. The degree of $u$, denoted by $deg(u)$,  is defined to be   the total number  of all occurrences of all $x\in X$ and $\theta\in \Omega$ in $u$.  For any associative  $\Omega$-word $u$, define   the depth of $u$ to be
$$
dep(u):=\min\{n|u\in \langle \Omega;  X \rangle_n \}.
$$
Let $(u)$ be a non-associative $\Omega$-word. Define the depth of $(u)$ by $dep((u))=dep(u)$.

For example, if
$
u=\omega^{(3)}(x_2x_1x_1, x_1, \omega^{(1)} (x_2  x_2x_1))x_2x_1,
$ where $x_1, x_2\in X$ and $\omega^{(3)}, \omega^{(1)} \in \Omega$,
  then $deg(u)=11$, $bre(u)=3$ and $dep(u)=2$.

Let $
u=u_1u_2\cdots u_n, n\geq 1,
$
where each $u_i$ is prime. Denote
$$
wt(u):=(deg(u), bre(u), u_1,u_2,\cdots, u_n).
$$
Let $X$ and $\Omega$ be well-ordered sets with the orders $>_X$ and $>_{\Omega}$, respectively.  Define the  Deg-lex order  $ >_{_{Dl}} $ on $\langle \Omega;  X\rangle$ as follows. For any $u=u_1u_2\cdots u_n, v=v_1v_2\cdots v_m$,   where    $u_i,  v_j$  are  prime, define
$$
u>_{_{Dl}}v \ \mbox{if}\ wt(u)>wt(v)\  \mbox{lexicographically},
$$
where if $u_i=\omega (u_{i1}, u_{i2}, \cdots, u_{it}), v_i=\theta(v_{i1}, v_{i2}, \cdots, v_{il})$ and $deg(u_i)=deg(v_i)$, then $u_i>v_i$  if
$$
 (\omega , u_{i1}, u_{i2}, \cdots, u_{it} )>(\theta, v_{i1}, v_{i2}, \cdots, v_{il})\  \mbox{lexicographically}.
$$

Let $\succ$ be the restriction of $>_{_{Dl}}$  on $X \cup \Omega( \langle \Omega;  X \rangle)$.
We define the     Lyndon-Shirshov  $\Omega$-words on the set $X$ by induction on the depth of the $\Omega$-words.

   For $n=0$, define
$$
ALSW(\Omega;  X)_0:=ALSW(X),
$$
$$
NLSW(\Omega;  X)_0:=NLSW(X)= \{[u]|  u\in ALSW(\Omega;  X)_0\}
$$
with respect to the lex-order $\succ_{lex}$ on $X^*$.

Assume that we have defined
$$
ALSW(\Omega;  X)_{n-1}\ \mbox{ and }\
NLSW(\Omega;  X)_{n-1}= \{[ u ]|u\in ALSW(\Omega;  X)_{n-1}\}.
$$

Define
$$
ALSW(\Omega;  X)_n:=ALSW(X\cup \Omega(ALSW(\Omega;  X)_{n-1}))
$$
with respect to the lex-order $\succ_{lex}$.

For  any    $u\in X\cup \Omega(ALSW(\Omega;  X)_{n-1}) $, define the bracketing way on $u$ by
$$
[u]:= \left\{
 \begin{array}{ll}
u,  & if \  u\in X, \\
\omega^{(m)}([ u_1], [u_2], \cdots,[ u_m]), & if \  u= \omega^{(m)}(u_1, u_2, \cdots, u_m).
 \end{array}
  \right.
$$
Denote
$$
[ X\cup \Omega(ALSW(\Omega;  X)_{n-1})]:=\{[ u]|u\in X\cup \Omega(ALSW(\Omega;  X)_{n-1})\}.
$$
Then, the order $\succ$ on $X\cup \Omega(ALSW(\Omega;  X)_{n-1})$ induces an order (still denoted  by $\succ$) on  $[ X\cup \Omega(ALSW(\Omega;  X)_{n-1})]$   by $[ u] \succ[ v ]$  if $u\succ v$ for any $u, v\in X\cup \Omega(ALSW(\Omega;  X)_{n-1})$.  If $u=u_1u_2\cdots u_m\in ALSW(\Omega;  X)_n$, where each  $u_i\in X\cup \Omega(ALSW(\Omega;  X)_{n-1})$, then we define
$$
[u]: =[[ u_1 ][ u_2]\cdots [ u_m ]]
$$
is the Shirshov standard  bracketing way  on  $\{[ u_1 ],[ u_2],\cdots ,[ u_m ]\}$, which means $[u]$ a non-associative  Lyndon-Shirshov   word on  $\{[ u_1 ],[ u_2],\cdots ,[ u_m ]\}$   with respect to the lex-order $\succ_{lex}$.
Denote
\begin{eqnarray*}
NLSW(\Omega;  X)_n&:=&\{[  u ]|u\in ALSW(\Omega;  X)_n\}.
\end{eqnarray*}
It is easy to see that 
$$
NLSW(\Omega;  X)_n=NLSW([ X\cup \Omega(ALSW(\Omega;  X)_{n-1})]).
$$
Denote 
\begin{eqnarray*}
ALSW(\Omega;  X)&:=&\bigcup_{n=0}^{\infty}ALSW(\Omega;  X)_n,\\
  NLSW(\Omega;  X)&:=&\bigcup_{n=0}^{\infty} NLSW(\Omega;  X)_n.
\end{eqnarray*}
Then
$$
  NLSW(\Omega;  X) =\{[  u ]|u\in  ALSW(\Omega;  X)\}.
$$
The elements of $ALSW(\Omega;  X)$ (resp. $NLSW(\Omega;  X)$) are called the   associative (resp. non-associative) Lyndon-Shirshov  $\Omega$-words.
By the above definitions, for any associative  Lyndon-Shirshov  $\Omega$-word $u$, $u$ is associated a unique non-associative  Lyndon-Shirshov  $\Omega$-word $[u] $.
For example, if
$$
u=\omega^{(3)}(x_2x_1x_1, x_1, \omega^{(1)} (x_2  x_2x_1))x_2x_1,
$$ where $x_1, x_2\in X$ with  $x_2 > x_1$, and $\omega^{(3)}, \omega^{(1)} \in \Omega$,
then $u \in ALSW(\Omega;  X)$ and
  $$
  [u]= (\omega^{(3)}(((x_2x_1)x_1), x_1, \omega^{(1)} ((x_2  (x_2x_1))))(x_2x_1))\in  NLSW(\Omega;  X).
  $$

\subsection{Free   Lie $\Omega$-algebras}

In this subsection, we prove that the  set $NLSW(\Omega;  X)$  of  all non-associative  Lyndon-Shirshov $\Omega$-words is  a linear basis of the  free Lie $\Omega$-algebra on the set $X$.

An associative (resp. non-associtive,   Lie) $\Omega$-algebra  over a field $k$   is an associative (resp.  non-associtive,  Lie) algebra $A$ with   multilinear operators $\Omega$.

For any associative  $\Omega$-algebra
$(A, \cdot, \Omega)$,   it is easy to see that  $(A, [, ], \Omega)$ is a  Lie $\Omega$-algebra, where
$$
[a, a']=a\cdot a'-a'\cdot a,\   \ a, a'\in A.
$$

 Let $X$ be a set.
An associative (resp. non-associtive,   Lie) $\Omega$-algebra $F(X)$ together with a injective map $i: X\rightarrow F(X)$  is called a free   associative (resp. non-associtive,   Lie) $\Omega$-algebra on  $X$, if for any   associative (resp. non-associtive,   Lie) $\Omega$-algebra $A$ and any map $\sigma: X\rightarrow A$, there exists a unique associative (resp. non-associtive,   Lie) $\Omega$-algebra homomorphism $\tilde{\sigma}: F(X)\rightarrow A$  such that $\tilde{\sigma}i=\sigma$.

Let  $k\langle \Omega;  X\rangle$  be the  $k$-linear space spanned by  $\langle \Omega;  X\rangle$.  Then $k\langle \Omega;  X\rangle$ is a  free  associative $\Omega$-algebra on the set $X$, see  \cite{bcq}. Denote  $Lie(\Omega;  X)$   the   Lie  $\Omega$-subalgebra of  $k\langle \Omega;  X\rangle$  generated by $X$ under the Lie bracket
$
(uv)=uv-vu.
$
The elements of $k\langle \Omega;  X\rangle$ (resp. $Lie(\Omega;  X)$) are called associative (resp. Lie)  $\Omega$-polynomials   on $X$.

\begin{lemma}\label{le2.7}
For any  $(u)\in  (\Omega;  X)$,   $(u)$ has a representation
$$
(u)=\sum \alpha_i[u_i],
$$
in $Lie(\Omega;  X)$, where each $\alpha_i\in k$ and $[u_i]\in NLSW(\Omega;  X)$.
\end{lemma}
{\bf Proof.} Induction on $dep((u))$. If $dep((u))=0$, by Lemma \ref{le2.1}, we have
$
(u)=\sum \alpha_i[u_i],
$
where each $\alpha_i\in k$ and $[u_i]\in NLSW(X)\subseteq NLSW(\Omega;  X)$.

Assume that the result is true for any $(u)$ with  $dep((u))\leq n-1$.

Let $dep((u))= n\geq 1$. There are two cases to consider.

Case 1. If $(u)=\theta((u_1),(u_2),\cdots,(u_m))$, then $dep((u_i))\leq n-1, 1 \leq i \leq m$. Thus, by induction,  we may assume that
$
(u_i)= [u_i],
$
where each  $[u_i]\in NLSW(\Omega;  X)$.
Therefore,
$$
(u)=\theta([u_1],[u_2],\cdots,[u_m]) = [\theta(u_1,u_2,\cdots,u_m)] \in NLSW(\Omega;  X).
$$

Case 2. If
$
(u)=(a_1a_2\cdots a_m), m\geq 2,
$
where each  $a_i$ is prime,   then by   Case 1, we may assume that
$
(u)=([a_1][a_2]\cdots [a_m]),
$
where each $a_i \in X\cup \Omega(ALSW(\Omega;  X)_{n-1})$. Thus, $(u)$ is a non-associative word on $\{ [a_1],[a_2],\cdots, [a_m]\}$.
By   Lemma \ref{le2.1}, we have
$$
([a_1][a_2]\cdots [a_m])=\sum \alpha_i[u_i],
$$
where  each  $\alpha_i\in k$  and  $[u_i]\in  NLSW( \{ [a_1],[a_2],\cdots, [a_m]\})\subseteq NLSW(\Omega;  X)$.  \hfill $ \square$\\

For any $f\in k\langle \Omega;  X\rangle$, let $\overline{f}$ be the   leading $\Omega$-word  of $f$ with respect to the order $>_{_{Dl}}$ on $\langle \Omega;  X\rangle$.

\begin{lemma}\label{lem2.8}
If  $u\in ALSW(\Omega;  X)$, then $\overline{[u]}=u$ with respect to  the order $>_{_{Dl}}$ on $\langle \Omega;  X\rangle$.
\end{lemma}
{\bf Proof.} Induction on $dep(u)$. If $dep(u)=0$, then $\overline{[u]}=u$ by Lemma \ref{le2.2}. Assume that the result is true for any $u\in ALSW(\Omega;  X)$ with  $dep(u)\leq n-1$.

Let $dep(u)= n\geq 1$. There are two cases to consider.

Case 1. If $bre(u)=1$  and $u= \theta(u_1,u_2,\cdots,u_m)$, then   $[u]=\theta([u_1],[u_2],\cdots,[u_m])$. By induction,  we have
$
\overline{[u_i]}  =u_i, 1\leq i\leq m.
$
Thus,
$$
\overline{[u]}=\theta(\overline{[u_1]},\overline{[u_2]},\cdots,\overline{[u_m]})= \theta(u_1,u_2,\cdots,u_m)=u.
$$

Case 2. If $bre(u)>1$  and
$
u=a_1a_2\cdots a_m,
$
where each $a_i$ is prime,  then by Lemma \ref{le2.2}, we have
$$
[u]=[[a_1][a_2]\cdots [a_m]] =[a_1][a_2]\cdots [a_m]+\sum \alpha_i [a_{i_1}][a_{i_2}]\cdots [a_{i_m}],
$$
where each $\{ i_1, i_2, \ldots, i_m\}=\{1,2, \ldots, m\}$ and   $a_{i_1}a_{i_2}\cdots a_{i_m} <_{_{Dl}}  a_1a_2\cdots a_m$. By Case 1, we have $\overline{[a_t]}=a_t $ and $\overline{[a_{i_j}]}=a_{i_j}$. It follows that
$
\overline{[u]}=a_1a_2\cdots a_m=u.
$
\hfill $ \square$\\

\begin{lemma}\label{le2.9}
$NLSW(\Omega;  X)$ is a linear basis of $Lie(\Omega;  X)$.
\end{lemma}
{\bf Proof.} Suppose
$
\sum_{i=1}^{m}\alpha_i[u_i]=0
$
in $Lie(\Omega;  X)$, where each $\alpha_i\in k$,  $u_i\in ALSW(\Omega;  X)$ and $u_i >_{_{Dl}} u_{i+1}$. If $\alpha_1\neq 0$, then by Lemma \ref{lem2.8}, we have
$
\overline{\sum_{i=1}^{m}\alpha_i[u_i]}=u_1,
$
a contradiction. Therefore, $NLSW(\Omega;  X)$ is  linear independent set. By Lemma \ref{le2.7},  $NLSW(\Omega;  X)$ is a linear basis of $Lie(\Omega;  X)$.\hfill $ \square$\\

Let  $k(\Omega;  X)$  be the $k$-linear space spanned by   $(\Omega;  X)$. It is easy to see that  $k(\Omega;  X)$ is a free non-associtive  $\Omega$-algebra on the set $X$.
Denote $R_\Omega$  the set consisting of the following relations in $k(\Omega;  X)$:
$$
((u)(v))=-((v)(u)),
$$
$$
(((u)(v))(w))=(((u)(w))(v))+( (u)((v)(w))),
$$
where $(u), (v), (w)\in (\Omega;  X)$.
Then
$
FL_\Omega(X):=k(\Omega;  X)/Id(R_\Omega)
$
is a free  Lie $\Omega$-algebra on the set $X$.

\begin{theorem}
$Lie(\Omega;  X)$ is a free  Lie $\Omega$-algebra on the set $X$ with   a linear basis $NLSW(\Omega;  X)$.
\end{theorem}
{\bf Proof.}
Let $i:X\rightarrow FL_\Omega(X), x\mapsto x+Id(R_\Omega)$ and $\varphi: X\rightarrow k\langle \Omega;  X\rangle, x\mapsto x$. Since $ FL_\Omega(X)$ is a free Lie $\Omega$-algebra on $X$, there is a unique Lie $\Omega$-algebra homomorphism $\widetilde{\varphi}:  FL_\Omega(X)\rightarrow  k\langle \Omega;  X\rangle$ such that $\widetilde{\varphi}i=\varphi$. It is easy to see that
$$
\widetilde{\varphi}(FL_\Omega(X))=Lie(\Omega;  X).
$$
Similar to  proof of Lemma \ref{le2.7}, we have for any $(u)\in (\Omega;  X)$,
$$
(u)+Id(R_\Omega)=\sum \beta_j[u_j]+Id(R_\Omega),
$$
where each $\beta_j\in k$ and  $u_j\in ALSW(\Omega;  X)$.
If
$$
\widetilde{\varphi}(\sum_{i=1}^{m}\alpha_i[u_i]+Id(R_\Omega))=\sum_{i=1}^{m}\alpha_i[u_i]=0
$$
in $k\langle \Omega;  X\rangle$,
 where each $\alpha_i\in k$,  $u_i\in ALSW(\Omega;  X)$, then by Lemma \ref{le2.9}, we have each  $\alpha_i=0$. Thus, $\widetilde{\varphi}$ is injective. It follows that $FL_\Omega(X)\cong Lie(\Omega;  X)$, i.e.,
 $Lie(\Omega;  X)$ is a free  Lie $\Omega$-algebra on the set $X$.
 \hfill $ \square$\\

\section{Gr\"{o}bner-Shirshov bases   for   Lie $\Omega$-algebras}

\subsection{Composition-Diamond lemma for  associative $\Omega$-algebras}

In this subsection, we review    Gr\"{o}bner-Shirshov
bases theory   for    associative  $\Omega$-algebras, which can be found in  \cite{bcq}.

Let $k\langle \Omega;  X\rangle$ be the free   associative $\Omega$-algebra on $X$ and $\star\notin X$. By a
$\star$-$\Omega$-word we mean any expression in $\langle \Omega; X\cup
\{\star\}\rangle$ with only one occurrence of $\star$. The set of all
  $\star$-$\Omega$-word on $X$ is denoted by $\langle \Omega;  X \rangle^\star$.

Let $\pi$ be a   $\star$-$\Omega$-word and $s\in k\langle
\Omega;  X\rangle$. Then we call
$$
\pi|_{s}:=\pi|_{\star\mapsto s}
$$
an   $s$-$\Omega$-word.

Now, we assume that $\langle \Omega;  X\rangle$ is equipped with a monomial
order  $>$. This means that $>$ is a well order  on
$\langle \Omega;  X\rangle$ such that for any $ v, w \in \langle \Omega;  X\rangle$ and
$\pi\in \langle \Omega;  X\rangle^\star$, if $w> v$, then  $\pi|_w> \pi|_v$.

For every    $\Omega$-polynomial $f\in k\langle \Omega;  X\rangle $, let
$\bar{f}$  be the leading   $\Omega$-word of $f$ with respect to the order $> $. If the coefficient
of $\bar{f}$ is $1$, then we call that $f$ is  monic. We also call  a  set $S \subseteq k\langle \Omega;  X\rangle $   monic if each  $s\in S$ is monic.

\ \

Let $f, g\in k\langle \Omega;  X\rangle $ be   monic. Then we define   two
kinds of compositions.
\begin{enumerate}
\item[(I)]If there exists an associative $\Omega$-word $w=\bar{f}a=b\bar{g}$ for some $a,b\in
\langle \Omega;  X\rangle$ such that $bre(w)< bre(\bar{f})+bre(\bar{g})$, then
we call $(f,g)_{w}:=fa-bg$ the intersection composition of $f$
and $g$ with respect to the ambiguity $w$.
\item[(II)] If there exists an  associative $\Omega$-word  $w=\bar{f}=\pi|_{\bar{g}}$ for some
$\pi \in \langle \Omega;  X\rangle^\star$, then we call $(f,g)_{w}:=f-\pi|_{g}$ the
inclusion  composition of  $f$ and $g$ with respect to the ambiguity $w$.
\end{enumerate}

Let $S \subseteq k\langle \Omega;  X\rangle$ be   monic. The composition
$(f,g)_w$ is called trivial modulo $(S,w)$ if
$$
(f,g)_w=\sum\alpha_i\pi_i|_{s_i},
$$
where each $\alpha_i\in k$,  $\pi_i\in \langle \Omega;  X\rangle^\star$, $s_i\in
S$ and $\pi_i|_{\overline{s_i}}< w$. If this is the case, we write
$$
(f,g)_w\equiv_{ass} 0 \ mod (S,w).
$$

In general, for any two associative $\Omega$-polynomials $p$ and $q$, $ p\equiv_{ass}
q \ mod (S,w) $ means that $ p-q=\sum\alpha_i\pi_i|_{s_i}, $ where
each $\alpha_i\in k$,  $\pi_i\in \langle \Omega;  X\rangle^\star$, $s_i\in S$
and $\pi_i|_{\overline{s_i}}< w$.

A monic set $S$ is called a Gr\"{o}bner-Shirshov basis  in  $k\langle
\Omega;  X\rangle$ if any composition $(f,g)_w$ of $f,g\in S$  is
trivial modulo $(S,w)$.

\begin{lemma}\label{th3.1} {\em(\cite{bcq}, Composition-Diamond lemma  for   associative  $\Omega$-algebras)}\ \  Let $S \subseteq k\langle \Omega;  X\rangle$ be   monic,   $> $ a monomial order  on $\langle \Omega;  X\rangle$  and  $Id_{ass}(S)$ the ideal of $k\langle \Omega;  X\rangle$ generated by $S$.   Then the following
statements are equivalent:
 \begin{enumerate}
\item[(i)] $S $ is a Gr\"{o}bner-Shirshov basis in $k\langle \Omega;  X\rangle$.
\item[(ii)] $ f\in Id_{ass}(S)\Rightarrow \bar{f}=\pi|_{\overline{s}}$
for some $\pi \in \langle \Omega;  X\rangle^\star$ and $s\in S$.
\item[(iii)]The set $Irr(S) = \{ w\in \langle \Omega;  X\rangle |  w \neq
\pi|_{\overline{s}},
  \pi \in \langle \Omega;  X\rangle^\star,    s\in S\}$
is a  linear basis of  the associative $\Omega$-algebra $k\langle \Omega;  X|S\rangle:=k\langle
\Omega;  X\rangle/Id_{ass}(S)$.
\end{enumerate}
\end{lemma}

\subsection{Composition-Diamond lemma    for   Lie  $\Omega$-algebras}

In this subsection, we  establish   Composition-Diamond lemma for   Lie  $\Omega$-algebras, which is a generalization of the Shirshov's Composition-Diamond lemma for   Lie algebras.

Let $>_{_{Dl}}$ be the order defined as before. It is easy to see that  $>_{_{Dl}}$ is a monomial order on $\langle \Omega;  X\rangle$. We always use this order in this subsection, in particular, a Gr\"{o}bner-Shirshov basis  in  $k\langle
\Omega;  X\rangle$ is with respect to the order $>_{_{Dl}}$.

The following lemma follows from Lemma \ref{lem2.3}.
\begin{lemma} Let $u=u_1u_2\cdots u_m$, where each $u_i\in X \cup \Omega(ALSW(\Omega;  X))$.  Then there  exists  a unique decomposition
$$
u=c_1c_2\cdots c_t,
$$
where each $c_i\in ALSW(\Omega;  X)$ and $c_j\preceq_{lex}c_{j+1}$, $1 \leq j\leq t-1$.
\end{lemma}

\begin{lemma} Let  $\pi\in \langle \Omega;  X\rangle^\star$ and $v, \pi|_v\in ALSW(\Omega;  X)$.  Then there  is a  $\pi_{_1}\in \langle \Omega;  X\rangle^\star$ and $c\in \langle \Omega;  X\rangle$ such that
$$
[\pi|_v]=[\pi_{_1}|_{[vc]}],
$$
where $c$ may be empty. Let
$$
[\pi|_v]_{v}:=[\pi_{_1}|_{[vc]}]|_{[vc]\mapsto [\cdots[[[v][c_1]][c_2]]\cdots [c_m]]},
$$
where   $c=c_1c_2\cdots c_m$ with  each  $c_i\in ALSW(\Omega;  X)$ and  $c_t \preceq_{lex} c_{t+1}, 1\leq t\leq m-1$.  Then,
$$
[\pi|_v]_{v}=\pi|_{[v]}+\sum \alpha_i \pi_i|_{[v]},
$$
where each $\alpha_i\in k$ and $\pi_i|_{v}<_{_{Dl}}\pi|_{v}$. It follows that $ \overline{[\pi|_v]_{v}}=\pi|_v$   with respect to the order $>_{_{Dl}}$.
\end{lemma}
{\bf Proof.} Induction on the depth of $\pi|_v$.  If $dep(\pi|_v)=0$, then the result  is true by Lemma \ref{le2.5}.
Assume that the result is true for any $\pi|_v$ with $dep(\pi|_v)\leq n-1$.

Let  $dep(\pi|_v)=n\geq 1$. There are two cases to consider.

Case 1. If $\pi|_v=avb$,  where $a, b\in \langle \Omega;  X\rangle$. Let $a=a_1a_2\cdots a_l, l\geq 0$ and $b=b_1b_2\cdots b_t, t\geq 0$, where   $  a_i$ and $b_j$ are prime.
By  Lemma \ref{le2.5}, we can obtain
$$
[\pi|_v]=[[a_1][a_2]\cdots [a_l][vc] [b_{j+1}]\cdots [b_t] ] =[a[vc]d],
$$
where $c=b_1b_2\cdots b_j$ and $d=b_{j+1}b_{j+2}\cdots b_t$.
Let  $c=c_1c_2\cdots c_m$, where each  $c_i\in ALSW(\Omega;  X)$ and  $  c_1 \preceq_{lex} \cdots \preceq_{lex} c_m$.
Then,
$$
[\pi|_v]_v=[[a_1][a_2]\cdots [a_l][\cdots[[[v][c_1]][c_2]]\cdots[c_m]]] [b_{j+1}]\cdots [b_t] ]
$$
and by   Lemma \ref{le2.5},  we have
$$
[\pi|_v]_v=[a_1][a_2]\cdots [a_l][v] [b_1][b_2]\cdots [b_t] +\sum \beta_i [d_{i_1}] \cdots [d_{i_p}][v] [d_{i_{p+1}}] \cdots [d_{i_{p+t}}].
$$
Therefore,  by Lemma  \ref{lem2.8},  we have
$$
[\pi|_v]_v=a[v]b+\sum \alpha_j a_j'[v]b_j',
$$
where each $ a_j'vb_j'<_{_{Dl}}avb$.

Case 2. If $\pi|_v=a \theta(u_1, u_2, \cdots, \pi'|_{v}, \cdots, u_q)b$, where $a, b\in \langle \Omega;  X\rangle$ and may be empty, then   we have
$$
[\pi|_v]= [a \theta(u_1, u_2, \cdots, [\pi'|_{v}], \cdots, u_q)b].
$$
By induction,
$$
[\pi'|_{v}]=[\pi''|_{[vc]}], \ \
[\pi'|_{v}]_v=\pi'|_{[v]}+\sum \alpha_{i}\pi_i'|_{[v]},
$$
where $\pi_i'|_{v} <_{_{Dl}} \pi'|_{v} $. Therefore,
\begin{eqnarray*}
[\pi|_v]&=& [a \theta(u_1, u_2, \cdots, [\pi''|_{[vc]}], \cdots, u_q)b]
\end{eqnarray*}
and
\begin{eqnarray*}
 {[\pi|_v]_{v}} &=&[a \theta(u_1,   \cdots, \pi'|_{v}, \cdots, u_q)b]_v \\
&=&[a  \theta([u_1],   \cdots, [\pi'|_{v}]_v, \cdots, [u_q]) b]\\
&=&a \theta([u_1],   \cdots, [\pi'|_{v}]_v, \cdots, [u_q])b+\sum \alpha_i a_i \theta([u_1],  \cdots, [\pi'|_{v}]_v, \cdots, [u_q])b_i\\
&=& a\theta(u_1,   \cdots, \pi'|_{[v]}, \cdots, u_q)b+ \sum \beta_j a_j'\theta(u_{j_1}',  u_{j_2}' \cdots, \pi''_j|_{[v]}, \cdots, u_{j_q}')b_j'\\
&=& \pi|_{[v]}+\sum \beta_j a_j'\theta(u_{j_1}', u_{j_2}',  \cdots, \pi''_j|_{[v]}, \cdots, u_{j_q}')b_j',
\end{eqnarray*}
where each $ a_j'\theta(u_{j_1}',   \cdots, \pi''_j|_{v}, \cdots, u_{j_q}')b_j'<_{_{Dl}} \pi|_v$. It follows that
$\overline{[\pi|_v]_v}=\pi|_v$.\hfill $ \square$\\

\begin{definition}
Let $\pi\in \langle \Omega;  X\rangle^\star$ and $f\in Lie(\Omega;  X)\subseteq k\langle \Omega;  X\rangle$ be  monic. If $\pi|_{\bar{f}} \in ALSW(\Omega;  X)$, then
$$
[\pi|_{f}]_{_{\bar{f}}}:=[\pi|_{\bar{f}}]_{_{\bar{f}}}|_{_{[\bar{f}]\mapsto f}}
$$
is called a special normal $f$-word.
\end{definition}

\begin{corollary}\label{co2.14}
Let $f\in Lie(\Omega;  X)$ and $\pi|_{\bar{f}}\in ALSW(\Omega;  X)$. Then
$$
[\pi|_{f}]_{_{\bar{f}}}=\pi|_{f}+\sum\alpha_i \pi_i|_{f},
$$
\end{corollary}
where each $\alpha_i\in k$ and $\pi_i|_{_{\bar{f}}} <_{_{Dl}}  \pi|_{_{\bar{f}}}$.\\

Let $f, g \in Lie(\Omega;  X)$ be  monic.  There are two kinds of compositions.
\begin{enumerate}
\item[(a)]\
If  $w=\bar{f}a=b\bar{g}$, where $a,b\in \langle \Omega;  X\rangle$ with  $bre(w)<bre(\bar{f})+bre(\bar{g})$, then
$$
\langle f,g\rangle_w:=[fa]_{_{\overline{f}}} - [bg]_{_{\overline{g}}}
$$
is called the intersection composition of $f$ and $g$ with respect to the ambiguity $w$.

\item[(b)]\  If $w=\overline{f}=\pi|_{\overline{g}}$,  then
$$
\langle f,g\rangle_w:=f-[\pi|_{g}]_{_{\overline{g}}}
$$
is called the inclusion composition of $f$ and $g$ with respect to  the ambiguity  $w$.
\end{enumerate}

If $S$ is a monic subset of $Lie(\Omega;  X)$, then
the composition $\langle f,g\rangle_w$ is called trivial modulo $(S, w)$ if
$$
\langle f,g\rangle_w=\sum \alpha_i [\pi_i|_{s_i}]_{_{\overline{s_i}}},
$$
where each $\alpha_i\in k,
  \ s_i\in S$, $[\pi_i|_{s_i}]_{_{\overline{s_i}}}$  is a special  normal $s_i$-word and $\pi_i|_{_{\overline{s_i}}}<_{_{Dl}} w$. If
this is the case, then we write
$$
\langle f,g\rangle_w\equiv 0\  mod(S,w).
$$
In general, for any two Lie $\Omega$-polynomials $p$ and $q$, $ p\equiv
q \ mod (S,w) $ means that $ p-q=\sum\alpha_i[\pi_i|_{s_i}]_{_{\overline{s_i}}}, $ where
each $\alpha_i\in k,
  \ s_i\in S$, $[\pi_i|_{s_i}]_{_{\overline{s_i}}}$  is a special  normal $s_i$-word and $\pi_i|_{_{\overline{s_i}}}<_{_{Dl}} w$.

\begin{definition}
A monic set $S$ is called a
Gr\"{o}bner-Shirshov basis  in   $Lie(\Omega;  X)$ if any
composition $\langle f,g\rangle_w$ of  $f, g\in S$ is trivial modulo $(S,w)$.
\end{definition}

\begin{lemma}\label{l2.13}
 Let $f,g\in Lie(\Omega;  X)\subset{k\langle \Omega;  X\rangle}$ be monic. Then
$$
\langle f,g\rangle_w-(f,g)_w\equiv_{ass}0 \ \ mod(\{f,g\},w).
$$
\end{lemma}
{\bf Proof.} If $\langle f,g\rangle_w$ and $(f,g)_w$ are compositions of intersection, where $w=\overline{f}a=b\overline{g}$, then
$$
\langle f,g\rangle_w=[fa]_{\overline{f}} - [bg]_{\overline{g}}=fa+\sum \alpha_i a_ifa_i' -bg-\sum \beta_jb_jgb_j',
$$
where $a_i\bar{f}a_i', b_j\bar{g}b_j'<_{_{Dl}}  w$.  It follows that
$$
\langle f,g\rangle_w-(f,g)_w\equiv _{ass} 0 \ mod(\{f,g\}, w).
$$

If $\langle f,g\rangle_w$ and $(f,g)_w$ are compositions of inclusion, where $w=\bar{f}=\pi|_{\bar{g}}$, then
$$
\langle f,g\rangle_w=f-[\pi|_{g}]_{\bar{g}}=f-\pi|_{g}-\sum \alpha_i\pi_i|_{g},
$$
where $\pi_i|_{\bar{g}}<_{_{Dl}} w$.
It follows that
$$
\langle f,g\rangle_w-(f,g)_w\equiv _{ass} 0 \ mod(\{f,g\}, w). \ \ \ \ \ \ \ \  \hfill  \square
$$

\begin{lemma}\label{th3.5}
Let $S\subset{Lie(\Omega;  X)}\subset{k\langle \Omega;  X\rangle}$ be monic.  With the order $>_{_{Dl}}$ on $\langle \Omega; X\rangle$,    the following two  statements are equivalent:
\begin{enumerate}
\item[(i)]
$S$ is a Gr\"{o}bner-Shirshov basis
in $Lie(\Omega;  X)$.
\item[(ii)] $S$ is a Gr\"{o}bner-Shirshov basis in
$k\langle \Omega;  X\rangle$.
\end{enumerate}
\end{lemma}
{\bf Proof.} $(i)\Rightarrow (ii)$.
Suppose that $S$ is a Gr\"{o}bner-Shirshov basis in $Lie(\Omega;  X)$. Then,
for any composition $\langle f,g\rangle_w$, we have
$
\langle
f,g\rangle_w=\sum\alpha_i [\pi_i|_{s_i}]_{_{\overline{s_i}}},
$
where each $\alpha_i \in k, \pi_i|_{_{\overline{s_i}}}<_{_{Dl}} w, \pi_i\in \langle \Omega;  X\rangle^\star, s_i\in S$. By Corollary \ref{co2.14}, we have
$
\langle f,g\rangle_w=\sum\beta_t\pi_t|_{s_t},
$
where each $\beta_t\in k, \pi_t|_{_{\overline{s_t}}}<_{_{Dl}} w$. Therefore,  by Lemma \ref{l2.13}, we  can obtain that
$
(f,g)_w\equiv_{ass}0 \ mod(S,w).
$
Thus, $S$ is a Gr\"{o}bner-Shirshov basis in $k\langle \Omega;  X\rangle$.

$(ii)\Rightarrow(i)$. Assume that $S$ is a Gr\"{o}bner-Shirshov basis in
$k\langle \Omega;  X\rangle$. Then, for any composition $\langle
f,g\rangle_w$ in $S$,    we  have  $\langle
f,g\rangle_w\in Lie(\Omega;  X)$ and $\langle
f,g\rangle_w\in Id_{ass}(S)$.  By Composition-Diamond lemma for associative $\Omega$-algebras,   $\overline{\langle
f,g\rangle_w}=\pi_1|_{_{\overline{s_1}}}\in ALSW(\Omega;  X)$.
Let
$$
h_1=\langle
f,g\rangle_w-\alpha_1[\pi_1|_{s_1}]_{_{\overline{s_1}}},
$$
 where $\alpha_1$ is the coefficient of $\overline{\langle
f,g\rangle_w}$.
Then,
$\overline{h_1}<_{_{Dl}} \overline{\langle f,g\rangle_w}$, $h_1\in Id_{ass}(S)$ and $h_1\in Lie(\Omega;  X)$. Now, the result follows   from   induction on $\overline{\langle f,g\rangle_w}$.   \hfill $\square$\\

\begin{lemma}\label{le3.6}
Let  $S\subset{Lie(\Omega;  X)}$ be  monic and
$$
Irr(S):= \{[w]|w\in ALSW(\Omega;  X), w\neq \pi|_{\bar{s}},  s\in S,  \pi\in\langle \Omega;  X\rangle^\star\}.
$$
Then, for any  $h\in{Lie(\Omega;  X)}$, $h$ can be expressed as
$$
h=\sum\alpha_i[u_i]+
\sum\beta_j[\pi_j|_{s_j}]_{_{\overline{s_j}}},
$$
where each $u_i\in ALSW(\Omega;  X), u_i\leq_{_{Dl}} \bar{h}$ and $s_i\in S$, $\pi_j|_{_{\overline{s_j}}}\leq_{_{Dl}} \bar{h}$.
\end{lemma}
{\bf Proof.} Since $f\in Lie(\Omega;  X)$,
$
h=\sum  \alpha_{i}[u_{i}],
$
where each $u_{i}\in ALSW(\Omega;  X)$ and  $u_{i}>_{_{Dl}} u_{i+1}$.
If $[u_1]\in Irr(S)$, then let
$
h_{1}=h-\alpha_{1}[u_1].
$
Otherwise,   there exists $s_1\in{S}$ such that $u_1=\pi|_{_{\overline{s_1}}}$. Let
$
h_1=h-[\pi|_{s_1}]_{_{\overline{s_1}}}.
$
In both of  the above two cases, we have  $h_1\in{Lie(\Omega;  X)}$ and $h_1<_{_{Dl}} \bar{h}$. Then, by induction on $\bar{h}$,  we can obtain the result.
 \hfill $ \square$\\

 The following theorem is   Composition-Diamond lemma  for   Lie $\Omega$-algebras. It is an analogue  of  Shirshov's Composition
lemma for Lie algebras \cite{Sh}, which was specialized to
associative algebras by L. A. Bokut \cite{Bo76}, see also G.M.  Bergman
\cite{Be78} and B. Buchberger \cite{bu65, Bu70}.

\begin{theorem}\label{cdll}
(Composition-Diamond lemma for   Lie $\Omega$-algebras) Let $S\subseteq Lie(\Omega;  X)$ be a non-empty monic set and $Id_{Lie}(S)$ the ideal of $Lie(\Omega;  X)$ generated by $S$.
 Then the following statements are equivalent:
\begin{enumerate}
\item[(i)] $S $ is a Gr\"{o}bner-Shirshov basis  in $Lie(\Omega;  X)$.
\item[(ii)] $f\in Id_{Lie}(S)\Rightarrow \bar{f}=\pi|_{\bar{s}} \in ALSW(\Omega;  X)$ for some $s\in S$ and $\pi\in\langle \Omega;  X\rangle^\star$.
\item[(iii)] The set
$$
Irr(S)= \{[w]|w\in ALSW(\Omega;  X), w\neq \pi|_{\bar{s}}, s\in S,  \pi\in\langle \Omega;  X\rangle^\star \}
$$
is a  linear basis of  the  Lie $\Omega$-algebra $Lie(\Omega;  X|S):=Lie(\Omega;  X)/Id_{Lie}(S)$.
\end{enumerate}
\end{theorem}
{\bf Proof.} $(i) \Rightarrow(ii)$ \  Since $f\in Id_{Lie}(S)\subseteq Id_{ass} (S)$, by Lemmas  \ref{th3.5} and   \ref{th3.1}, we have
$\bar{f}=\pi|_{\bar{s}}\in ALSW(\Omega;  X)$ for some $s\in S$ and $\pi\in\langle \Omega;  X\rangle^\star$.

$(ii)\Rightarrow(iii)$   \ Suppose that
$
\sum\alpha_i[u_i]=0
$
in $Lie(\Omega;  X|S)$, where each $[u_i]\in Irr(S)$ and  $u_i>_{_{Dl}}  u_{i+1}$.  That
is,
$
\sum\alpha_i[u_i]\in{Id_{Lie}(S)}.
$
Then each $\alpha_i$ must be 0. Otherwise, say $\alpha_1\neq0$,
since
$
\overline{\sum\alpha_i[u_i]}=u_1
$
and by (ii), we have  $u_1\notin Irr(S)$,    a contradiction. Therefore, $Irr(S)$ is linear independent. By Lemma \ref{le3.6},   $Irr(S)$
is a  linear basis of     $Lie(\Omega;  X|S)$.

$(iii)\Rightarrow(i)$  \ For any composition $\langle
f,g\rangle_w$  of $f,g\in{S}$, we have $\langle
f,g\rangle_w\in{Id_{Lie}(S)}$. Then, by (iii) and   Lemma \ref{le3.6},
$$
 \langle
f,g\rangle_w=\sum\beta_j[\pi_j|_{s_j}]_{_{\overline{s_j}}},
$$
where each $\beta_j\in k,   \pi_j\in\langle \Omega;  X\rangle^\star, s_j\in S, \pi_j|_{_{\overline{s_j}}} <_{_{Dl}} w$. This proves that $S$ is a
Gr\"{o}bner-Shirshov basis in $Lie(\Omega;  X)$.   \hfill $ \square$\\

\section{Applications}
In this  section, as applications of Theorem \ref{cdll}, we give Gr\"{o}bner-Shirshov bases for   free $\lambda$-Rota-Baxter Lie algebras, free    modified $\lambda$-Rota-Baxter Lie algebras and free Nijenhuis Lie algebras and then  linear bases of such three free  algebras  are obtained.

\subsection{Free     $\lambda$-Rota-Baxter  Lie algebras}

A Lie algebra $L$ equiped with a  linear map $ P:L\rightarrow L$ satisfying
$$
 [P(x)P(y)]  = P(  [P(x) y] )  +P( [x P(y) ])+\lambda P ( [x y ]), \  x, y \in  L
$$
 is called a     Rota-Baxter  Lie algebra of weight $\lambda$ or  $\lambda$-Rota-Baxter  Lie algebra.
 It is easy to see that any  $\lambda$-Rota-Baxter  Lie algebra $(L, P)$ is a Lie $\Omega$-algebra, where $\Omega=\{P\}$.

Let $Lie(\{P\}; X)$   be the free   Lie   $ \{P\}$-algebra   on   a set $X$.  Denote $S$   the  set consisting of the following  Lie   $\{P\}$-polynomials in $Lie(\{P\}; X)$:
$$
f_{u, v} =(P([u])P([v]))- P((P([u])[v])) - P(([u]P([v])))-\lambda P(([u][v])),
$$
where $u, v \in ALSW(\{P\}; X)$ and $u>_{_{Dl}}v$.

It is easy to see
$$
RBL(X):=Lie(\{P\}; X|S)=Lie(\{P\}; X)/Id_{Lie}(S)
$$
is the free $\lambda$-Rota-Baxter Lie algebra on the  set $X$.

\begin{theorem}\label{th4.1}With the order  $>_{_{Dl}}$  on $\langle \{P\}; X\rangle$,  the set $S$ is a Gr\"{o}bner-Shirshov basis in $Lie(\{P\}; X)$. It follows that the set
$$
Irr(S)=  \left\{
[w]\in NLSW(\{P\}; X)
 \left|
 \begin{array}{ll}
  w\neq \pi|_{P(u)P(v)}, \ \    \pi\in   \langle \{P\}; X\rangle^\star \\
 u, v \in ALSW(\{P\}; X), u >_{_{Dl}} v  \\
\end{array}
\right. \right\}
$$
is a linear basis of the free   $\lambda$-Rota-Baxter  Lie algebra $RBL(X)=Lie(\{P\}; X|S)$ on   $X$.
\end{theorem}
{\bf  Proof.} All the  possible compositions of
Lie  $ \{P\}$-polynomials in $S$ are listed as below:
$$\langle f_{u,v}, f_{v, w}\rangle_{w_1}, w_1=P(u)P(v)P(w), u >_{_{Dl}}v >_{_{Dl}} w,$$
$$\langle f_{\pi|_{P(u)P(v)},w}, f_{u, v}\rangle_{w_2}, w_2=P(\pi|_{P(u)P(v)})P(w),   u>_{_{Dl}} v,  \pi|_{P(u)P(v)} >_{_{Dl}} w,$$
$$\langle f_{u, \pi|_{P(v)P(w)}}, f_{v, w}\rangle_{w_3}, w_3=P(u)P(\pi|_{P(v)P(w)}),  v >_{_{Dl}}  w,  u>_{_{Dl}}\pi|_{P(v)P(w)}.$$

We check that all the compositions are trivial.
\begin{eqnarray*}
&&\langle f_{u,v}, f_{v, w}\rangle_{w_1}\\
&=&[f_{u,v}P(w)]_{_{\overline{f_{u,v}}}}-[P(u)f_{v,w}]_{_{\overline{f_{v,w}}}} \\
&=& ((P([u])P([w]))P([v]))-  (P((P([u])[v]))P([w])) -  (P(([u]P([v])))P([w])) \\
&&-  \lambda (P(([u][v]))P([w]))+ (P([u]) P((P([v])[w]))) +  (P([u]) P(([v]P([w]))))\\
&&+\lambda (P([u]) P(([v][w]))).
\end{eqnarray*}

By direct computation, we have  $mod(S, w_1)$:
\begin{eqnarray*}
& & ((P([u])P([w]))P([v]))\\
&\equiv&( P((P([u])[w]))+ P(([u]P([w])))+\lambda P(([u][w])))P([v]))\\
&\equiv& P(( P((P([u])[w]))[v]))+P(((P([u])[w])P([v])))+\lambda P((P([u])[w]))[v]))\\
&&+\  P((P(([u]P([w])))[v]))+P((([u]P([w]))P([v])))+\lambda P((([u]P([w]))[v])) \\
&&+ \lambda P((P(([u][w]))[v]))+\lambda P((([u][w])P([v])))+\lambda^2P((([u][w])[v])),
\end{eqnarray*}
\begin{eqnarray*}
&& (P(([u]P([v])))P([w]))\\
&\equiv& P((P(([u]P([v])))[w]))+P((([u]P([v]))P([w])))+\lambda P((([u]P([v]))[w]))\\
&\equiv& P((P(([u]P([v])))[w]))+P((([u]P([w]))P([v])))+P(([u](P([v])P([w])))\\
&&+\lambda P((([u]P([v]))[w]))\\
&\equiv& P((P(([u]P([v])))[w]))+P((([u]P([w]))P([v])))+P(([u]P((P([v])[w]))))\\
&&+P(([u]P(([v] P([w])))))+\lambda P(([u]P(([v] [w]))))
+\lambda P((([u]P([v]))[w])),
\end{eqnarray*}
\begin{eqnarray*}
&&(P((P([u])[v]))P([w]))\\
&\equiv& P((P((P([u])[v]))[w]))+P(((P([u])[v])P([w])))+\lambda P(((P([u])[v])[w]))\\
&\equiv&  P((P((P([u])[v]))[w]))+P(((P([u])P([w]))[v]))+P((P([u])([v]P([w]))))\\
&&+\lambda P(((P([u])[v])[w]))\\
&\equiv& P((P((P([u])[v]))[w]))+P((P((P([u])[w]))[v]))+P((P(([u]P([w])))[v]))\\
&&+\lambda P((P(([u][w]))[v]))+P((P([u])([v]P([w]))))+\lambda P(((P([u])[v])[w])),
\end{eqnarray*}
\begin{eqnarray*}
(P(([u][v]))P([w]))
\equiv P((P(([u][v]))[w]))+P((([u][v])P([w])))+\lambda P((([u][v])[w])),
\end{eqnarray*}
\begin{eqnarray*}
&&(P([u])P((P([v])[w])))\\
&\equiv& P((P([u])(P([v])[w])))+P(([u]P((P([v])[w])))+\lambda P(([u](P([v])[w])))\\
&\equiv& P((P([u])P([v]))[w]))+P((P([v])(P([u])[w])))
 + P(([u](P((P([v])[w])))\\
  &&+ \lambda P(([u](P([v])[w])))\\
&\equiv& P((P((P([u])[v]))[w])) + P((P(([u]P([v])))[w]))+\lambda  P((P(([u][v]))[w]))\\
&&+P((P([v])(P([u])[w])))+P(([u](P((P([v])[w]))))) +\lambda P(([u](P([v])[w]))),
\end{eqnarray*}
\begin{eqnarray*}
&&(P([u])P(([v]P([w])))\\
&\equiv& P((P([u])([v]P([w]))))+P(([u]P(([v]P([w])))))+\lambda P(([u]([v]P([w]))),
\end{eqnarray*}
\begin{eqnarray*}
( P([u])P(([v] [w])))
\equiv  P((P([u])([v] [w])))+P(([u]P(([v] [w]))))+\lambda P(([u]([v][w]))).
\end{eqnarray*}

Therefore, we have
$$
\langle f_{u,v}, f_{v, w}\rangle_{w_1}=[f_{u,v}P(w)]_{_{\overline{f_{u,v}}}}-[P(u)f_{v,w}]_{_{\overline{f_{v,w}}}}\equiv 0\ mod(S, w_1).
$$

Denote
$$
r(\pi|_{f_{u, v}})=[\pi|_{f_{u, v}}]_{_{\overline{f_{u,v}}}}-[\pi|_{_{\overline{f_{u,v}}}}].
$$
Then,
\begin{eqnarray*}
&&\langle f_{\pi|_{P(u)P(v)},w}, f_{u, v}\rangle_{w_2}\\
&=& f_{\pi|_{P(u)P(v)},w}-[P(\pi|_{f_{u, v}})P(w)]_{P(u)P(v)} \\
&=&f_{\pi|_{P(u)P(v)},w}-(P([\pi|_{f_{u, v}}]_{P(u)P(v)})P(w))\\
&=& P((P([\pi|_{P(u)P(v)}])[w]))+P(([\pi|_{P(u)P(v)}]P([w])))+\lambda P(([\pi|_{P(u)P(v)}][w]))\\
&&- (P(r(\pi|_{f_{u, v}}))P([w]))\\
&\equiv& P((P(r(\pi|_{f_{u, v}}))[w]))+P((r(\pi|_{f_{u, v}})P([w])))
 + \lambda P((r(\pi|_{f_{u, v}})[w]))\\
&&-P((P(r(\pi|_{f_{u, v}}))[w]))-P(( r(\pi|_{f_{u, v}})  P([w])))
-\lambda P((r(\pi|_{f_{u, v}}) [w]))\\
&\equiv& 0\ mod(S, w_2),
\end{eqnarray*}
 \begin{eqnarray*}
&&\langle f_{u, \pi|_{P(v)P(w)}}, f_{v, w}\rangle_{w_3}\\
&=& f_{u, \pi|_{P(v)P(w)}}-[P(u)P(\pi|_{f_{v, w}})]_{P(v)P(w)} \\
&=&f_{u, \pi|_{P(v)P(w)}}-(P(u)P([\pi|_{f_{v, w}}]_{P(v)P(w)}))\\
&=& P((P([u]) [\pi|_{P(v)P(w)}]))+P(([u]P([\pi|_{P(v)P(w)}] )))+\lambda P(( [u][\pi|_{P(v)P(w)}]))\\
&&- (P([u])P(r(\pi|_{f_{v, w}})))\\
&\equiv&   P((P([u]) r(\pi|_{f_{v, w}})))+P(([u]P(r(\pi|_{f_{v, w}}))))+\lambda P(( [u]r(\pi|_{f_{v, w}})))\\
&&-P((P([u]) r(\pi|_{f_{v, w}})))-P(([u]P(r(\pi|_{f_{v, w}}))))-\lambda P(( [u]r(\pi|_{f_{v, w}})))\\
&\equiv& 0\ mod(S, w_3).
\end{eqnarray*}


Therefore,   $S$ is a Gr\"{o}bner-Shirshov basis in $Lie(\{P\}; X)$.
By Composition-Diamond lemma for Lie $\Omega$-algebras, the set $Irr(S)$
 is a linear basis of the free   $\lambda$-Rota-Baxter  Lie algebra $RBL(X)$.
\hfill $ \square$\\

\subsection{Free  modified $\lambda$-Rota-Baxter Lie algebras}

A Lie algebra $L$ equipped with a  linear map $ P:L\rightarrow L$ satisfying
$$
 [P(x)P(y)]  = P(  [P(x) y] )  +P( [x P(y) ])+\lambda  [x y ] , \  x, y \in  L
$$
is called a     modified $\lambda$-Rota-Baxter Lie algebra. If $\lambda=-1$, the  modified $\lambda$-Rota-Baxter Lie algebra is also called a  Baxter Lie algebra in  \cite{bor90}.
 It is easy to see that any modified $\lambda$-Rota-Baxter Lie algebra   $(L, P)$ is a Lie $\Omega$-algebra, where $\Omega=\{P\}$.

Let $Lie(\{P\}; X)$   be the free   Lie   $ \{P\}$-algebra   on   a set $X$. Denote $S_B$   the  set consisting of the following  Lie  $\{P\}$-polynomials in $Lie(\{P\}; X)$:
$$
f_{u, v}^B =(P([u])P([v]))- P((P([u])[v])) - P(([u]P([v])))-\lambda ([u][v]),
$$
where $u, v \in ALSW(\{P\}; X)$ and $u>_{_{Dl}}v$.

It is easy to see
$$
MRBL(X):=Lie(\{P\}; X|S_B)=Lie(\{P\}; X)/Id_{Lie}(S_B)
$$
is the free modified $\lambda$-Rota-Baxter Lie algebra on   $X$.

\begin{theorem}\label{th4.2}With the order  $>_{_{Dl}}$  on $\langle \{P\}; X\rangle$,  the set $S_B$ is a Gr\"{o}bner-Shirshov basis in $Lie(\{P\}; X)$. It follows that the set
$$
Irr(S_B)=  \left\{
[w]\in NLSW(\{P\}; X)
 \left|
 \begin{array}{ll}
  w\neq \pi|_{P(u)P(v)}, \ \    \pi\in   \langle \{P\}; X\rangle^\star \\
 u, v \in ALSW(\{P\}; X), u >_{_{Dl}} v  \\
\end{array}
\right. \right\}
$$
is a linear basis of the free  modified $\lambda$-Rota-Baxter Lie algebra $MRBL(X)=Lie(\{P\}; X|S_B)$ on  $X$.
\end{theorem}
{\bf  Proof.} All the  possible compositions of
Lie  $ \{P\}$-polynomials in $S$ are listed as below:
$$\langle f_{u,v}^B, f_{v, w}^B\rangle_{w_1}, w_1=P(u)P(v)P(w), u >_{_{Dl}}v >_{_{Dl}} w,$$
$$\langle f^B_{\pi|_{P(u)P(v)},w}, f_{u, v}^B\rangle_{w_2}, w_2=P(\pi|_{P(u)P(v)})P(w),   u>_{_{Dl}} v,  \pi|_{P(u)P(v)} >_{_{Dl}} w,$$
$$\langle f^B_{u, \pi|_{P(v)P(w)}}, f^B_{v, w}\rangle_{w_3}, w_3=P(u)P(\pi|_{P(v)P(w)}),  v >_{_{Dl}}  w,  u>_{_{Dl}}\pi|_{P(v)P(w)}.$$


We check that all the compositions are trivial.
\begin{eqnarray*}
&&\langle f^B_{u,v}, f^B_{v, w}\rangle_{w_1}\\
&=&[f^B_{u,v}P(w)]_{_{\overline{f^B_{u,v}}}}-[P(u)f^B_{v,w}]_{_{\overline{f^B_{v,w}}}} \\
&=&(((P([u])P([v]))- P((P([u])[v])) - P(([u]P([v])))-\lambda ([u][v]))P(w))\\
&& -(P(u)(((P([v])P([w]))- P((P([v])[w])) - P(([v]P([w])))-\lambda ([v][w]))))  \\
&=& ((P([u])P([w]))P([v]))-  (P((P([u])[v]))P([w])) -  (P(([u]P([v])))P([w])) \\
&&-  \lambda (([u][v])P([w]))+ (P([u])P((P([v])[w]))) +  (P([u]) P(([v]P([w]))))\\
&&+\lambda (P([u])([v][w])).
\end{eqnarray*}
By direct computation, we have  $mod(S, w_1)$:
\begin{eqnarray*}
& & ((P([u])P([w]))P([v]))\\
&\equiv&( (P((P([u])[w]))+ P(([u]P([w])))+\lambda   ([u][w]))  P([v]))\\
&\equiv& P(( P((P([u])[w]))[v]))+P(((P([u])[w])P([v])))+\lambda ((P([u])[w]) [v])\\
&&+\  P((P(([u]P([w])))[v]))+P((([u]P([w]))P([v])))+\lambda  (([u]P([w]))[v])\\
&&+ \lambda  (([u][w])  P([v])),
\end{eqnarray*}
\begin{eqnarray*}
&& (P(([u]P([v])))P([w]))\\
&\equiv& P((P(([u]P([v])))[w]))+P((([u]P([v]))P([w])))+\lambda  (([u]P([v]))[w]),
\end{eqnarray*}
\begin{eqnarray*}
&&(P((P([u])[v]))P([w]))\\
&\equiv& P((P((P([u])[v]))[w]))+P(((P([u])[v])P([w])))+\lambda  ((P([u])[v])[w]),
\end{eqnarray*}
\begin{eqnarray*}
&&(P([u])P((P([v])[w])))\\
&\equiv& P((P([u])(P([v])[w])))+P(([u]P((P([v])[w])))+\lambda  ([u](P([v])[w])),
\end{eqnarray*}
\begin{eqnarray*}
&&(P([u])P(([v]P([w])))\\
&\equiv& P((P([u])([v]P([w]))))+P(([u]P(([v]P([w])))))+\lambda  ([u]([v]P([w])),
\end{eqnarray*}
\begin{eqnarray*}
&&P((([u]P([v]))P([w])))-P((([u]P([w]))P([v])))=P(([u](P([v])P([w]))))\\
&\equiv& P(([u]P((P([v])[w])))) +  P(([u]P(([v]P([w])))))+\lambda P(([u] ([v][w]))),
 \end{eqnarray*}
\begin{eqnarray*}
&& P(((P([u])[v])P([w])))-P((P([u])([v]P([w]))))=P(((P([u])P([w]))[v]))\\
 &\equiv& P(( P((P([u])[w])) [v]))+P((P(([u]P([w])) )[v]))+\lambda P((([u][w])[v])),
\end{eqnarray*}
\begin{eqnarray*}
&&P(((P([u])[w])P([v])))-   P((P([u])(P([v])[w])))=P(((P([u])P([v]))[w]))\\
&\equiv&  P((P((P([u])[v]))[w]))+P((P(([u] P([v])))[w]))+\lambda P(([u][v])[w])).
\end{eqnarray*}

Therefore, we have
$$
\langle f^B_{u,v}, f^B_{v, w}\rangle_{w_1}=[f^B_{u,v}P(w)]_{_{\overline{f^B_{u,v}}}}-[P(u)f^B_{v,w}]_{_{\overline{f^B_{v,w}}}}\equiv 0\ mod(S, w_1).
$$
Denote
$$
r(\pi|_{f^B_{u, v}})=[\pi|_{f^B_{u, v}}]_{_{\overline{f^B_{u,v}}}}-[\pi|_{_{\overline{f^B_{u,v}}}}].
$$
Then,
\begin{eqnarray*}
&&\langle f^B_{\pi|_{P(u)P(v)},w}, f^B_{u, v}\rangle_{w_2}\\
&=& f^B_{\pi|_{P(u)P(v)},w}-[P(\pi|_{f^B_{u, v}})P(w)]_{P(u)P(v)} \\
&=&f^B_{\pi|_{P(u)P(v)},w}-(P([\pi|_{f^B_{u, v}}]_{P(u)P(v)})P(w))\\
&=& P((P([\pi|_{P(u)P(v)}])[w]))+P(([\pi|_{P(u)P(v)}]P([w])))+\lambda  ([\pi|_{P(u)P(v)}][w]) \\
&&- (P(r(\pi|_{f^B_{u, v}}))P([w]))\\
&\equiv& P((P(r(\pi|_{f^B_{u, v}}))[w]))+P((r(\pi|_{f^B_{u, v}})P([w])))
 + \lambda  (r(\pi|_{f^B_{u, v}})[w]) \\
&&-P((P(r(\pi|_{f^B_{u, v}}))[w]))-P(( r(\pi|_{f^B_{u, v}}) P([w])))
-\lambda   (r(\pi|_{f^B_{u, v}}) [w]) \\
&\equiv& 0\ mod(S, w_2),
\end{eqnarray*}
\begin{eqnarray*}
&&\langle f^B_{u, \pi|_{P(v)P(w)}}, f^B_{v, w}\rangle_{w_3}\\
&=& f^B_{u, \pi|_{P(v)P(w)}}-[P(u)P(\pi|_{f^B_{v, w}})]_{P(v)P(w)} \\
&=&f^B_{u, \pi|_{P(v)P(w)}}-(P(u)P([\pi|_{f^B_{v, w}}]_{P(v)P(w)}))\\
&=& P((P([u]) [\pi|_{P(v)P(w)}]))+P(([u]P([\pi|_{P(v)P(w)}] )))+\lambda  ( [u][\pi|_{P(v)P(w)}]) \\
&&- (P([u])P(r(\pi|_{f_{v, w}})))\\
&\equiv&   P((P([u]) r(\pi|_{f_{v, w}})))+P(([u]P(r(\pi|_{f_{v, w}}))))+\lambda  ( [u]r(\pi|_{f_{v, w}})) \\
&&-P((P([u]) r(\pi|_{f_{v, w}})))-P(([u]P(r(\pi|_{f_{v, w}}))))-\lambda  ( [u]r(\pi|_{f_{v, w}})) \\
&\equiv& 0\ mod(S, w_3).
\end{eqnarray*}

Therefore,   $S_B$ is a Gr\"{o}bner-Shirshov basis in $Lie(\{P\}; X)$.
By Composition-Diamond lemma for Lie $\Omega$-algebras, the set $Irr(S_B)$
 is a linear basis of the free modified $\lambda$-Rota-Baxter Lie algebra $MRBL(X)$.
\hfill $ \square$\\

 \subsection{Free  Nijenhuis Lie algebras}

A Lie algebra $L$ equiped with a  linear map $ P:L\rightarrow L$ satisfying
$$
 [P(x)P(y)]  = P(  [P(x) y] )  +P( [x P(y) ])-P^2([x y ]) , \  x, y \in  L
$$
is called a     Nijenhuis   Lie algebra.

Let $Lie(\{P\}; X)$   be the free   Lie   $ \{P\}$-algebra   on   a set $X$. Denote $S_N$   the  set consisting of the following Lie   $\{P\}$-polynomials in $Lie(\{P\}; X )$:
$$
f^N_{u, v} =(P([u])P([v]))- P((P([u])[v])) - P(([u]P([v])))+P^2(([u][v])),
$$
where $u, v \in ALSW(\{P\}; X)$ and $u>_{_{Dl}}v$.

It is easy to see that
$$N
L(X):=Lie(\{P\}; X|S_N)=Lie(\{P\}; X)/Id_{Lie}(S_N)
$$
is the free Nijenhuis  Lie algebra on   $X$.

\begin{theorem}\label{th4.2}With the order  $>_{_{Dl}}$  on $\langle \{P\}; X\rangle$,  the set $S_N$ is a Gr\"{o}bner-Shirshov basis in $Lie(\{P\}; X)$. It follows that the set
$$
Irr(S_N)=  \left\{
[w]\in NLSW(\{P\}; X)
 \left|
 \begin{array}{ll}
  w\neq \pi|_{P(u)P(v)}, \ \    \pi\in   \langle \{P\}; X \rangle^\star \\
 u, v \in ALSW(\{P\}; X), u >_{_{Dl}} v  \\
\end{array}
\right. \right\}
$$
is a linear basis of the free    Nijenhuis   Lie algebra $NL(X)=Lie(\{P\}; X|S_N)$ on  $X$.
\end{theorem}
{\bf  Proof.} All the  possible compositions of
Lie  $ \{P\}$-polynomials in $S$ are listed as below:
$$\langle f^N_{u,v}, f^N_{v, w}\rangle_{w_1}, w_1=P(u)P(v)P(w), u >_{_{Dl}}v >_{_{Dl}} w,$$
$$\langle f^N_{\pi|_{P(u)P(v)},w}, f_{u, v}^N\rangle_{w_2}, w_2=P(\pi|_{P(u)P(v)})P(w),   u>_{_{Dl}} v,  \pi|_{P(u)P(v)} >_{_{Dl}} w,$$
$$\langle f^N_{u, \pi|_{P(v)P(w)}}, f^N_{v, w}\rangle_{w_3}, w_3=P(u)P(\pi|_{P(v)P(w)}),  v >_{_{Dl}}  w,  u>_{_{Dl}}\pi|_{P(v)P(w)}.$$


We check that all the compositions are trivial.
\begin{eqnarray*}
&&\langle f^N_{u,v}, f^N_{v, w}\rangle_{w_1}\\
&=&[f^N_{u,v}P(w)]_{_{\overline{f^N_{u,v}}}}-[P(u)f^N_{v,w}]_{_{\overline{f^N_{v,w}}}} \\
&=& ((P([u])P([w]))P([v]))-  (P((P([u])[v]))P([w])) -  (P(([u]P([v])))P([w])) \\
&&+ (P^2(([u][v]))P([w]))+ (P([u])P((P([v])[w])))+   (P([u])P(([v]P([w]))))\\
&&-  (P([u])P^2(([v][w]))).
\end{eqnarray*}

By direct computation, we have  $mod(S, w_1)$:
\begin{eqnarray*}
& & ((P([u])P([w]))P([v]))\\
&\equiv&( P((P([u])[w]))P([v]))+ (P(([u]P([w])))P([v]))- (P^2(([u][w]))P([v]))\\
&\equiv& P(( P((P([u])[w]))[v]))+P(((P([u])[w])P([v])))- P^2(((P([u])[w]) [v]))\\
&&+\  P((P(([u]P([w])))[v]))+P((([u]P([w]))P([v])))- P^2(( [u]P([w]) [v])) \\
&&-   P((P^2(([u][w]))[v]))-  P((P(([u][w]))P([v])))+ P^2((P(([u][w]))[v]))\\
&\equiv& P(( P((P([u])[w]))[v]))+P(((P([u])[w])P([v])))- P^2(((P([u])[w]) [v] ))\\
&&+\  P((P(([u]P([w])))[v]))+P((([u]P([w]))P([v])))- P^2((([u]P([w]) [v])) \\
&&-   P((P^2(([u][w]))[v]))-  P^2((P(([u][w]))  [v] ))- P^2(( ([u][w]) P( [v] )))\\
&&+P^3((([u][w])[v]))+ P^2((P(([u][w]))[v])),
\end{eqnarray*}
\begin{eqnarray*}
&& (P(([u]P([v])))P([w]))\\
&\equiv& P((P(([u]P([v])))[w]))+P((([u]P([v]))P([w])))- P^2((([u]P([v]))[w]))\\
&\equiv& P((P(([u]P([v])))[w]))+P((([u]P([w]))P([v])))+P(([u](P([v])P([w])))\\
&&- P^2((([u]P([v]))[w]))\\
&\equiv& P((P(([u]P([v])))[w]))+P((([u]P([w]))P([v])))+P(([u]P((P([v])[w]))))\\
&&+P(([u]P(([v] P([w])))))- P(([u]P^2(([v] [w]))))
-P^2((([u]P([v]))[w])),
\end{eqnarray*}
\begin{eqnarray*}
&&(P((P([u])[v]))P([w]))\\
&\equiv& P((P((P([u])[v]))[w]))+P(((P([u])[v])P([w])))- P^2(((P([u])[v])[w]))\\
&\equiv&  P((P((P([u])[v]))[w]))+P(((P([u])P([w]))[v]))+P((P([u])([v]P([w]))))\\
&&- P^2(((P([u])[v])[w]))\\
&\equiv& P((P((P([u])[v]))[w]))+P((P((P([u])[w]))[v]))+P((P(([u]P([w])))[v]))\\
&&- P((P^2(([u][w]))[v]))+P((P([u])([v]P([w]))))- P^2(((P([u])[v])[w])),
\end{eqnarray*}
\begin{eqnarray*}
&&(P^2(([u][v]))P([w]))\\
&\equiv& P((P^2(([u][v]))[w]))+P((P(([u][v]))P([w])))- P^2((P(([u][v]))[w]))\\
&\equiv& P((P^2(([u][v]))[w]))+ P^2((P(([u][v]))[w])) )+  P^2( ([u][v]) P([w])) \\
&& - P^3(( ([u][v]) [w])) - P^2((P(([u][v]))[w])),
\end{eqnarray*}
\begin{eqnarray*}
&&(P([u])P((P([v])[w])))\\
&\equiv& P((P([u])(P([v])[w])))+P(([u]P((P([v])[w])))- P^2(([u](P([v])[w])))\\
&\equiv& P((P([u])P([v]))[w]))+P((P([v])(P([u])[w])))
 + P(([u](P((P([v])[w])))\\
  &&-P^2(([u](P([v])[w])))\\
&\equiv& P((P((P([u])[v]))[w])) + P((P(([u]P([v])))[w]))-  P ((P^2(([u][v]))[w]))\\
&&+P((P([v])(P([u])[w])))+P(([u](P((P([v])[w]))))) -P^2(([u](P([v])[w]))),
\end{eqnarray*}
\begin{eqnarray*}
&&(P([u])P(([v]P([w])))\\
&\equiv& P((P([u])([v]P([w]))))+P(([u]P(([v]P([w])))))- P^2(([u]([v]P([w]))),
\end{eqnarray*}
\begin{eqnarray*}
&&( P([u])P^2(([v] [w])))\\
&\equiv&  P((P([u])P(([v] [w]))))+P(([u]P^2(([v] [w]))))- P^2(([u]P(([v][w]))))\\
&\equiv&    P^2((P([u])([v] [w]) )) + P^2(([u]P(([v] [w])))) - P^3(([u]([v] [w]))) \\
&&+P(([u]P^2(([v] [w]))))- P^2(([u]P(([v][w])))).
\end{eqnarray*}

Therefore, we have
$$
\langle f^N_{u,v}, f^N_{v, w}\rangle_{w_1}=[f^N_{u,v}P(w)]_{_{\overline{f^N_{u,v}}}}-[P(u)f^N_{v,w}]_{_{\overline{f^N_{v,w}}}}\equiv 0\ mod(S, w_1).
$$

Denote
$$
r(\pi|_{f^N_{u, v}})=[\pi|_{f^N_{u, v}}]_{_{\overline{f^N_{u,v}}}}-[\pi|_{_{\overline{f^N_{u,v}}}}].
$$
Then,
\begin{eqnarray*}
&&\langle f^N_{\pi|_{P(u)P(v)},w}, f^N_{u, v}\rangle_{w_2}\\
&=& f^N_{\pi|_{P(u)P(v)},w}-[P(\pi|_{f^N_{u, v}})P(w)]_{P(u)P(v)} \\
&=&f^N_{\pi|_{P(u)P(v)},w}-(P([\pi|_{f^N_{u, v}}]_{P(u)P(v)})P(w))\\
&=& P((P([\pi|_{P(u)P(v)}])[w]))+P(([\pi|_{P(u)P(v)}]P([w])))- P^2(([\pi|_{P(u)P(v)}][w]))\\
&&- (P(r(\pi|_{f^N_{u, v}}))P([w]))\\
&\equiv& P((P(r(\pi|_{f^N_{u, v}}))[w]))+P(( r(\pi|_{f^N_{u, v}}) P([w])))
 -P^2(( r(\pi|_{f^N_{u, v}}) [w]))\\
&&-P((P(r(\pi|_{f^N_{u, v}}))[w]))-P((  r(\pi|_{f^N_{u, v}} )  P([w])))
+ P^2(( r(\pi|_{f^N_{u, v}} ) [w]))\\
&\equiv& 0\ mod(S, w_2),
\end{eqnarray*}
 \begin{eqnarray*}
&&\langle f^N_{u, \pi|_{P(v)P(w)}}, f^N_{v, w}\rangle_{w_3}\\
&=& f^N_{u, \pi|_{P(v)P(w)}}-[P(u)P(\pi|_{f^N_{v, w}})]_{P(v)P(w)} \\
&=&f^N_{u, \pi|_{P(v)P(w)}}-(P(u)P([\pi|_{f^N_{v, w}}]_{P(v)P(w)}))\\
&=& P((P([u]) [\pi|_{P(v)P(w)}]))+P(([u]P([\pi|_{P(v)P(w)}] )))-P^2(( [u][\pi|_{P(v)P(w)}]))\\
&&- (P([u])P(r(\pi|_{f^N_{v, w}})))\\
&\equiv&   P((P([u]) r(\pi|_{f^N_{v, w}})))+P(([u]P(r(\pi|_{f^N_{v, w}}))))-P^2(( [u]r(\pi|_{f^N_{v, w}})))\\
&&-P((P([u]) r(\pi|_{f^N_{v, w}})))-P(([u]P(r(\pi|_{f^N_{v, w}}))))+ P^2(( [u]r(\pi|_{f^N_{v, w}})))\\
&\equiv& 0\ mod(S, w_3).
\end{eqnarray*}


Therefore,   $S_N$ is a Gr\"{o}bner-Shirshov basis in $Lie(\{P\}; X)$.
By Composition-Diamond lemma for Lie $\Omega$-algebras, the set $Irr(S_N)$
 is a linear basis of the free  Nijenhuis  Lie algebra $NL(X)$. \hfill $ \square$\\

\noindent{\bf Acknowledgement}:
The authors would like to express their deepest gratitude to Professor L.A. Bokut
  for  his  kind guidance, useful discussions and enthusiastic encouragements.

\end{document}